\documentclass[titlepage,11pt]{article}
\usepackage{latexsym}
\usepackage{amsfonts}
\usepackage{amsmath}
\usepackage{mathtools}
\newcommand\blackslug{\hbox{\hskip 1pt \vrule width 4pt height 8pt depth 1.5pt
        \hskip 1pt}}
\newcommand\bbox{\hfill \quad \blackslug \bigbreak}

\def\LL{,\ldots,}
\newcommand{\vare}{\varepsilon}

\newcommand{\cupcup}{\cup \cdots\cup}

\title{Induced subgraph density. II. Sparse and dense sets in cographs}
\author{
Jacob Fox\thanks{Supported by a Packard Fellowship and by NSF grants DMS-1953990 and DMS-2154129.}\\
Stanford University,\\ Stanford, CA 94305, USA
\and
Tung Nguyen\thanks{Supported by AFOSR grants
A9550-19-1-0187 and FA9550-22-1-0234, and by NSF grants  DMS-1800053 and DMS-2154169.}\\
Princeton University,\\ Princeton, NJ 08544, USA
\and
Alex Scott\thanks{Supported by EPSRC grant EP/X013642/1}\\
University of Oxford, \\
Oxford, UK
\and
Paul Seymour\thanks{Supported by AFOSR grants
A9550-19-1-0187 and FA9550-22-1-0234, and by NSF grants  DMS-1800053 and DMS-2154169.}\\
Princeton University,\\ Princeton, NJ 08544, USA}

\date{July 4, 2022; revised \today}

\newtheorem{thm}{}[section]

\newcommand{\Proof}{\noindent{\bf Proof.}\ \ }

\begin{document}
\maketitle
\begin{abstract}
A well-known theorem of R\"odl says that
for every graph 
$H$, and every $\vare>0$, there exists $\delta>0$ such that 
if $G$ does not contain an induced copy of $H$, then there exists $X\subseteq V(G)$
with $|X|\ge \delta|G|$ such that one of $G[X],\overline{G}[X]$ has edge-density at most $\vare$. 
But how  does $\delta$ depend on $\epsilon$?
Fox and Sudakov conjectured 
that the dependence is at most polynomial: that for all $H$ there exists $c>0$ such that for all $\vare$
with $0<\vare\le 1/2$, R\"odl's theorem holds with $\delta=\vare^c$. This conjecture implies the
Erd\H{o}s-Hajnal conjecture, and until now it had not been verified for any non-trivial graphs $H$. Our first result 
shows that it is true when $H=P_4$. Indeed, in that case we can take $\delta=\vare$, and insist that 
one of $G[X],\overline{G}[X]$ has maximum degree at most $\vare^2|G|$).

Second, we will show that every graph $H$ that can be obtained by substitution from copies of $P_4$ satisfies the Fox-Sudakov
conjecture. To prove this, we need 
to work with 
a stronger property. 
Let us say $H$ is {\em viral} if there exists $c>0$ such that for all $\vare$ with 
$0<\vare\le 1/2$, if $G$ contains at most $\vare^c|G|^{|H|}$ copies of $H$ as induced subgraphs, then there exists 
$X\subseteq V(G)$
with $|X|\ge \vare^c|G|$ such that one of $G[X],\overline{G}[X]$ has edge-density at most $\vare$. We will show that $P_4$
is  viral, using a ``polynomial $P_4$-removal lemma'' of Alon and Fox. 
We will also show that the class of viral graphs is closed under 
vertex-substitution.

Finally, we give a different strengthening of R\"odl's theorem: we show that if $G$ does not contain an induced copy of $P_4$, 
then its vertices can be partitioned into at most $480\vare^{-4}$ subsets $X$ such that one of $G[X],\overline{G}[X]$ has maximum
degree at most $\vare|X|$.

\end{abstract}

\section{Introduction}
Some terminology and notation: $G[X]$ denotes the
induced subgraph with vertex set $X$
of a graph $G$; $|G|$ denotes the number of vertices of $G$; $\overline{G}$ is the complement graph of $G$; $P_4$ denotes the path with four vertices; 
a graph is {\em $H$-free} if it has no induced subgraph isomorphic to $H$; and a {\em cograph} is a $P_4$-free graph.
The {\em edge-density} of a graph $G$ is its number of edges divided by $\binom{|G|}{2}$.

A very useful theorem of R\"odl~\cite{rodl} says:
\begin{thm}\label{rodlthm}
For every graph $H$ and every $\vare>0$, there exists
$\delta>0$ such that for every $H$-free graph $G$, there exists $X\subseteq V(G)$
with $|X|\ge \delta|G|$ such that one of $G[X],\overline{G}[X]$ has edge-density at most $\vare$.
\end{thm}

How does $\delta$ depend on $\vare$, for a given graph $H$? Fox and Sudakov~\cite{foxsudakov} proposed the 
conjecture that the dependence is polynomial:
\begin{thm}\label{foxconj}
{\bf Conjecture (\cite{foxsudakov}, conjecture 7.1):} For every graph $H$ there exists $c>0$ such that for every $\vare$ with $0<\vare\le 1/2$ and every $H$-free graph $G$,
there exists $X\subseteq V(G)$
with $|X|\ge \vare^c|G|$ such that one of $G[X],\overline{G}[X]$ has edge-density at most $\vare$.
\end{thm}
This conjecture is very strong, and until now had not been verified for any nontrivial graphs $H$.
It was motivated by the Erd\H{o}s-Hajnal conjecture~\cite{EH0,EH}, which it implies, but which we do not discuss 
here.\footnote{Since this paper was submitted for publication, 
Buci\'c, Fox and Pham~\cite{bucic} have proved
that all graphs $H$ satisfying the
Erd\H{o}s-Hajnal conjecture also satisfy the Fox-Sudakov conjecture, and indeed are viral (defined later).}

We first prove that 
\ref{foxconj} holds in a particularly nice form when $H=P_4$. We will show:
\begin{thm}\label{P4thm}
For every $\vare\in [0,1]$ and every cograph $G$, there exists $X\subseteq V(G)$
with $|X|\ge \vare |G|$ such that one of $G[X],\overline{G}[X]$ has maximum degree at most $\vare^2|G|$ (and so at most $\vare|X|$).
\end{thm}

We need to define ``vertex-substitution'' before we go on. Let $H_1,H_2$ be graphs, let $v\in V(H_1)$,
and let $N$ be the set of all neighbours of $v$ in $H_1$. Let $H$ be obtained from the disjoint union of $H_1\setminus \{v\}$
and $H_2$ by making every vertex of $H_2$ adjacent to every vertex in $N$. Then $H$ is obtained by {\em substituting $H_2$
for the vertex $v$ of $H_1$}, and this operation is called {\em vertex-substitution}.

We would like to prove that more graphs than just $P_4$ satisfy \ref{foxconj}, and one natural way is via vertex-substitution (for example, Alon, Pach and Solymosi \cite{aps} showed that graphs satisfying the Erd\H os-Hajnal conjecture are closed under vertex-substitution).
We have not been able to show that the graphs that satisfy \ref{foxconj} are closed under vertex-substitution. But we
{\em have} been able to show that $P_4$ itself has an even stronger property than \ref{foxconj}, and graphs with this stronger property are closed under vertex-substitution. Consequently:
\begin{thm}\label{foxcograph}
All graphs that can be obtained by vertex-substitution starting from copies of $P_4$ and its subgraphs satisfy \ref{foxconj}.
\end{thm}

Let us say a {\em copy} of $H$ in $G$ is an isomorphism from $H$
to an induced subgraph of $G$.
There is a theorem of Nikiforov~\cite{nikiforov}, strengthening R\"odl's theorem:
\begin{thm}\label{nikiforov}
For every graph $H$ and all $\vare>0$, there exists $\delta>0$ such that for every graph $G$, if there are at most 
$\delta|G|^{|H|}$ copies of $H$ in $G$,
then there exists $X\subseteq V(G)$ with $|X|\ge \delta|G|$
such that one of $G[X], \overline{G}[X]$ has edge-density at most $\vare$.
\end{thm}
Again, one could ask how $\delta$ depends on $\vare$. Let us say that $H$ is  {\em viral} if 
there exists $d>0$ such that for every graph $G$ and every $\vare$ with $0<\vare\le 1/2$,
either 
\begin{itemize}
\item there are at least $\vare^d|G|^{|H|}$ copies of $H$ in $G$; or
\item there exists $X\subseteq V(G)$
with $|X|\ge \vare^d|G|$ such that one of $G[X],\overline{G}[X]$ has edge-density at most~$\vare$.
\end{itemize}

We will show, using a recent ``polynomial removal lemma'' for $P_4$, proved by Alon and Fox~\cite{alonfox}, that:
\begin{thm}\label{P4virus}
All graphs with at most four vertices are viral.
\end{thm}
We will also show:
\begin{thm}\label{gemsubst}
If $H_1,H_2$ are viral and $H$ is obtained by substituting $H_2$
for a vertex of $H_1$, then $H$ is viral.
\end{thm}
We deduce:
\begin{thm}\label{viruscograph}
All graphs that can be obtained by vertex-substitution starting from graphs with at most four vertices are viral.
\end{thm}

In the final section, we  will discuss a different strengthening of R\"odl's theorem, and prove: 
\begin{thm}\label{P4partn}
If $G$ is a cograph, then for every $\vare$ with $0<\vare\le 1$,
there is a partition of $V(G)$ into at most $480\vare^{-4}$ sets such that for each of them, say $X$,
one of $G[X],\overline{G}[X]$ has maximum degree at most $\vare|X|$.
\end{thm}

\section{Cographs have large dense or sparse sets}

In this section we prove \ref{P4thm}.  We will discuss how close it is to best possible in the next section.

Cographs are well understood. There is a theorem  discovered independently by several authors 
(see~\cite{corneil}), that:
\begin{thm}\label{cograph}
If $G$ is a cograph with $|G|\ge 2$ then one of $G, \overline{G}$ is disconnected.
\end{thm}
We will use \ref{cograph} to prove \ref{P4thm} by induction on $|G|$. Applying it directly 
does not seem to work, and to use induction we will use a strengthening of \ref{P4thm},
the following (\ref{P4thm} follows by setting $x=y=\vare$):

\begin{thm}\label{betterthm}
If $G$ is a cograph then, for  all $x,y\ge 0$ with $\min(x,y)\le 1$, either:
\begin{itemize}
\item there exists $X\subseteq V(G)$ with $|X|\ge x|G|$ such that $G[X]$ has maximum degree at most $xy|G|$; or
\item there exists  $Y\subseteq V(G)$ with $|Y|\ge y|G|$ such that $\overline{G}[Y]$ has maximum degree at most $xy|G|$.
\end{itemize}
\end{thm}
\Proof If $|G|\le 1$ the result is true, so we assume that $|G|\ge 2$ and the result holds for all cographs with fewer vertices, 
and for all choices of $x,y\ge 0$ with $\min(x,y)\le 1$. If $x> 1$, then $y\le 1$ and the second bullet holds choosing $Y\subseteq V(G)$ 
with $|Y|=\lceil y|G|\rceil$; so we may assume that $x\le 1$ and similarly $y\le 1$. 
By \ref{cograph}, taking complements if necessary, we may assume that $G$ is not connected; let $G_1,G_2$ be two non-null subgraphs of 
$G$, with union $G$ and with $V(G_1)\cap V(G_2)=\emptyset$. Now we are given $x,y\ge 0$ with $\min(x,y)\le 1$. 

For $i = 1,2$, let $y_i=y|G|/|G_i|$. If for some $i\in \{1,2\}$ there exists $Y_i\subseteq V(G_i)$ with $|Y_i|\ge y_i|G_i|=y|G|$ 
such that 
$\overline{G_i}[Y_i]$ has maximum degree at most $xy_i|G_i|=xy|G|$, then the second bullet holds. Hence we assume that for 
$i = 1,2$ there is no such $Y_i$. But $\min(x,y_i)\le 1$, so from the inductive hypothesis, for $i = 1,2$
there exists $X_i \subseteq V(G_i)$ with $|X_i|\ge x|G_i|$ such that $G[X_i]$ has maximum degree at most $xy_i|G_i|=xy|G|$.
Then $|X_1\cup X_2|\ge x|G|$ and $G[X_1\cup X_2]$ has maximum degree at most $xy|G|$, and the first bullet of the theorem holds.
This proves \ref{betterthm}.~\bbox

Here is a consequence, strengthening \ref{P4thm}:
\begin{thm}\label{product}
Let $G$ be a cograph, and let $0\le \vare\le 1$. Then there exists $X,Y\subseteq V(G)$, such that $G[X], \overline{G}[Y]$ both have maximum degree at most
$\vare|G|$, and with $|X|\cdot |Y|\ge \vare|G|^2$.
\end{thm}
\Proof Let $I$ be the set of $x\in [0,1]$ such that for some $X\subseteq V(G)$,  $|X|\ge x|G|$ and $G[X]$ has maximum degree at most
$\vare|G|$; and let $J$ be the set of $x\in [0,1]$ such that for some $Y\subseteq V(G)$, $x|Y|\ge \vare|G|$ and $\overline{G}[Y]$ has maximum degree at most
$\vare|G|$. By \ref{betterthm}, $I\cup J=[0,1]$. Since $I,J$ are nonempty closed sets (because $G$ is finite), it follows that $I\cap J\ne \emptyset$.
This proves \ref{product}.~\bbox

The form of \ref{betterthm} seems novel, and suggests that we ask which other graphs have the same property.
Let us say $G$ is {\em good} if for all $x,y$ with $0\le x,y\le 1$, either:
\begin{itemize}
\item there exists $X\subseteq V(G)$ with $|X|\ge x|G|$ such that $G[X]$ has maximum degree at most $xy|G|$; or
\item there exists  $Y\subseteq V(G)$ with $|Y|\ge y|G|$ such that $\overline{G}[Y]$ has maximum degree at most $xy|G|$.
\end{itemize}
Thus, complements of good graphs are good; \ref{betterthm} says that all cographs are good; and its proof shows that goodness 
is preserved under taking disjoint unions. Which other graphs are good? This is still open, but we can show (we omit the proofs):
\begin{itemize}
\item  all forests are
good; 
\item the bull is not good; 
\item a cycle of length at least five is good if and only if its length is a multiple of six; and
\item goodness is {\em not} preserved under vertex-substitution; indeed, 
substituting a two-vertex graph for a vertex of a good graph does not always preserve goodness.
\end{itemize}

\section{The tightness of \ref{P4thm}}

Let us say two disjoint subsets $A,B$ are {\em complete} to each other if every vertex in $A$ is adjacent to every vertex in $B$, and {\em anticomplete}
if there are no edges between $A,B$.

For $\vare\in [0,1]$, let $\delta_{\vare}$ be the supremum of all $\delta$ such that
for every cograph $G$, there exists $X\subseteq V(G)$ such that $|X|\ge \delta|G|$ and one of $G[X],\overline{G}[X]$ has maximum degree at most
$\vare\delta|G|$. The next result shows that \ref{P4thm} is almost tight.
\begin{thm}\label{bounds}
If $\vare \in [0,1)$, then 
$\vare\le \delta_{\vare}\le (\lceil \vare^{-1}\rceil-1)^{-1}$.
\end{thm}
\Proof
By \ref{P4thm}, $\vare\le \delta_{\vare}$. 
Let $m= \lceil \vare^{-1}\rceil-1$; that is, the largest integer strictly less than $1/\vare$.
Take an integer $n\ge  (1/m-\vare)^{-1}$, and let $G$ be the cograph
consisting of the disjoint union of $m$ complete graphs $C_1\LL C_m$, each with $n$ vertices.
We will show that if $X\subseteq V(G)$ and one of $G[X],\overline{G}[X]$ has maximum degree at most $\vare |X|$, then
$|X|\le |G|/m$, and consequently $\delta_{\vare}\le 1/m$. Let $X\subseteq V(G)$ with $|X|>|G|/m=n$. By the pigeonhole principle, $|X\cap C_i|\ge |X|/m$ 
for some $i$, and so 
$G[X]$ has maximum degree at least $|X|/m - 1> \vare|X|$ (because $(1/m-\vare)|X|>   (1/m-\vare)n\ge 1$).
But since $|X|>|G|/m= |C_i|$, there is a vertex in $X\setminus C_i$, and the degree of
this vertex in $\overline{G}[X]$ is at least $|X|/m>\vare|X|$. This proves \ref{bounds}.~\bbox

The result \ref{P4thm} is neat, and one might think it should be tight, but it is not; and indeed, neither of the 
bounds of \ref{bounds} is tight when $1/2\le \vare<1$.
We will show that $\delta_\vare=1/(2-\vare)$ in this range. To do so, we first show
 the following, which implies that $\delta_\vare\ge 1/(2-\vare)>\vare$ when $1/2\le \vare<1$:
\begin{thm}\label{toprange}
Let $1/2\le \vare< 1$ and let $\delta=1/(2-\vare)$. For every non-null cograph $G$, there is a set $X\subseteq V(G)$
with $|X|>\delta |G|$ such that one of $G[X], \overline{G}[X]$ has maximum degree at most 
$\vare\delta|G|$. 
\end{thm}
\Proof
Let $G$ be a non-null cograph, and let $1/2\le \vare< 1$. Let $\delta=1/(2-\vare)$ and $d=\vare/(2-\vare)$; we must show that
there is a set $X\subseteq V(G)$
with $|X|>  \delta|G|$ such that one of $G[X], \overline{G}[X]$ has maximum degree 
at most $d|G|$.

We partition $V(G)$ into sets $X_1\LL X_k$ as follows. Suppose that $i\ge 1$ and we have 
defined $X_1\LL X_{i-1}$, such that $V(G)\ne X_1\cupcup X_{i-1}$. Let 
$Y=V(G)\setminus (X_1\cupcup X_{i-1})$. If 
$|Y|=1$, let $X_i=Y$ and $k=i$.
Now we assume that $|Y|>1$, and define $X_i$ as follows.
By \ref{cograph}, one of $G[Y],\overline{G}[Y]$ is not connected. Let $X_i$ be a
subset of $Y$ that is the vertex set of a component of one of $G[Y],\overline{G}[Y]$, chosen with 
$|X_i|$ minimum. 
Thus $|X_i|\le |Y|/2$, and in particular $V(G)\ne X_1\cupcup X_{i}$. This completes the inductive 
definition.
\\
\\
(1) {\em We may assume that $|X_i|\le \delta|G|/2$ for $1\le i\le k-1$.}
\\
\\
Suppose that some $|X_i|> \delta|G|/2$, and let $Y=V(G)\setminus  (X_1\cupcup X_{i})$. Choose $A\subseteq X_i$ with $|A|=\lfloor \delta|G|/2 +1\rfloor$.
As we saw, $|Y|\ge |X_i|$, and so there exists $B\subseteq Y$ with $|B|=|A|$. Now the set $A\cup B$ has cardinality more than $\delta|G|$. Moreover,
from the construction, $X_i$ is either complete or anticomplete to $Y$, and by
taking complements if necessary, we may assume the former. But then every vertex in $A$ has no neighbours in $B$ and has at most 
$|A|-1\le \delta|G|/2\le \vare\delta|G|$ neighbours in $A$, and similarly for $B$; and so setting $X=A\cup B$ satisfies the theorem. This proves (1).

\bigskip

We may assume that $|G|\ge 2$ and so $k\ge 2$. If $d|G|\ge |G|-1$, then the theorem is satisfied with $X=V(G)$ (because $\delta<1$ and every vertex 
has at most $d|G|$ neighbours in $G$). So we may assume that $d|G|< |G|-1$.
Choose $h$ with $0\le h\le k-1$, minimum such that $|X_{h+1}\cupcup X_k|\le d|G|+1$. (This is possible since the condition is satisfied when $h=k-1$).
Since $|G|>d|G|+1$ it follows that $h\ge 1$. By moving to the complement if necessary, we may assume that there $X_h, Y$ are anticomplete,
where $Y=X_{h+1}\cupcup X_k$. Let $I$ be the set of all $i\in \{1\LL h\}$ such that $X_i, Y$ are anticomplete, and let $J$
be the set of all $i\in \{1\LL h\}$ such that  $X_i, Y$ are complete. Thus $h\in I$. Moreover, all the sets $X_i\;(i\in I)$ are pairwise anticomplete,
and the sets $X_i\;(i\in J)$ are pairwise complete.

Choose $Z\subseteq X_h$ such that $|Y\cup Z|=\lfloor d|G|+1\rfloor$ (this is possible since $|X_h\cup Y|> d|G|+1$ from the minimality of $h$).
Let $A$ be the union of $Y$ and the sets $X_i\;(i\in I)$. 
Since each of the sets $X_i\; (i\in I)$ and $Y$ have cardinality at most $d|G|+1$ by (1), and there are no edges between them, it follows
that $G[A]$ has maximum degree at most $d|G|$. Similarly, let $B$ be the union of $Y\cup Z$ and the sets $X_i\;(i\in J)$; then since these sets all
have cardinality at most $d|G|+1$, and there are no edges
of $\overline {G}$ between any two of them, it follows that $\overline{G}[B]$ has maximum degree at most $d|G|$. But $|A|+|B|=|G|+|Y|+|Z|$,
and so one of $|A|,|B|$ has cardinality at least $(|G|+|Y|+|Z|)/2$. To complete the proof it suffices to show that $(|G|+|Y|+|Z|)/2\ge \delta|G|$.
Certainly $|Y\cup Z|> d|G|$; 
and hence 
$$(|G|+|Y|+|Z|)/2> (1+d)|G|/2=(1+\vare/(2-\vare))|G|/2=\delta|G|.$$
This proves \ref{toprange}.~\bbox

\ref{toprange} says that $|X|>\delta|G|$, and hence $|X|\ge \lfloor \delta|G|+1\rfloor$. Next we show that this is tight.
\begin{thm}\label{tightex}
Let $1/2\le \vare\le 1$ and $\delta=1/(2-\vare)$. For each even integer $2n\ge 4$, there is a cograph $G$ with $2n$ vertices such that 
if $X\subseteq V(G)$
and one of $G[X],\overline{G}[X]$ has maximum degree at most $\vare\delta|G|$, then $|X|\le \delta|G|+1$ (and hence $|X|\le \lfloor \delta|G|+1\rfloor$).
\end{thm}
\Proof
Let $G$ be the ``half-graph'' with vertex set $\{a_1\LL a_n, b_1\LL b_n\}$, in which $\{a_1\LL a_n\}$ is a stable set,
$\{b_1\LL b_n\}$ is a clique, and
$a_i, b_j$ are adjacent if and only if $i\le j$. This graph is a cograph. Now choose $X\subseteq V(G)$
such that $G[X]$ has maximum degree at most $\vare\delta|G|$, with $|X|$ maximum, and
subject to that with $|X\cap A|$ maximum. Since
$|X|>|G|/2$ (because $|G|\ge 4$ and so $\vare\delta|G|\ge 1$), $X$ contains a vertex $b\in B$. For each $a\in A\setminus X$, since $b$ is adjacent 
to all neighbours of $a$ in $X$, it follows that $(X\cup \{a\})\setminus \{b\}$ would be a better choice than $X$, a contradiction; and so 
$A\subseteq X$. Let $|X\cap B| =i$ say; then there is a vertex in $X\cap B$ with $i$ neighbours in $A$ and adjacent to all other vertices in
$X\cap B$, and since its degree in $G[X]$ is at most $\vare\delta|G|$, we deduce that $2i-1\le \vare\delta|G|$. So $|X\cap B|\le (\vare\delta|G|+1)/2$,
and hence $|X|\le |G|/2+ (\vare\delta|G|+1)/2=\delta|G|+1/2$. Similarly (the graph is not quite self-complementary), if $X\subseteq V(G)$
and $\overline{G}[X]$ has maximum degree at most $\vare\delta|G|$, it follows that $|X|\le \delta|G|+1$. This proves \ref{tightex}.~\bbox

We deduce:
\begin{thm}\label{exact}
If $1/2\le \vare\le 1$, then $\delta_{\vare} = 1/(2-\vare)$.
\end{thm}
\Proof
By \ref{toprange}, $\delta_{\vare}\ge 1/(2-\vare)$. By \ref{tightex}, $\delta_{\vare}\le 1/(2-\vare) + 1/(2n)$ for each integer $n\ge 2$,
and so 
$\delta_{\vare}\le 1/(2-\vare)$. This proves \ref{exact}.~\bbox

\section{Viral graphs and vertex-substitution}\label{sec:subst}
Let us 
prove \ref{gemsubst}, which we restate:
\begin{thm}\label{gemsubst2}
If $H_1,H_2$ are viral and $H$ is obtained by  substituting $H_2$
for a vertex of $H_1$, then $H$ is viral.
\end{thm}
\Proof
Let $H$ be obtained by  substituting $H_2$
for a vertex $v$ say of $H_1$. For $i = 1,2$, since $H_i$ is viral, 
there exists $d_i$ as in the definition of ``viral''. Let $d=(|H_2|+1)(d_1+1)+d_2$.
To show that $H$ is viral, we will show that: 
\\
\\
(1) {\em For every graph $G$ and all $\vare$ with $0<\vare\le 1/2$, either
\begin{itemize}
\item there exists $X\subseteq V(G)$
with $|X|\ge \vare^d|G|$ such that one of $G[X],\overline{G}[X]$ has edge-density at most $\vare$; or
\item there are at least $\vare^d|G|^{|H|}$ copies of $H$ in $G$.
\end{itemize}
}
\noindent
(We remind the reader that ``copy'' means an isomorphism from $H$ to an induced subgraph of $G$.)
Since $\vare^{d_1}|G|\ge \vare^{d}|G|$, we may assume that there is no $X\subseteq V(G)$
with $|X|\ge \vare^{d_1}|G|$ such that one of $G[X],\overline{G}[X]$ has edge-density at most $\vare$, since otherwise
the first bullet of (1) holds. Consequently,
from the choice of $d_1$, 
there are at least $\vare^{d_1}|G|^{|H_1|}$ copies of $H_1$ in $G$.
For each copy $\phi$ of $H_1\setminus \{v\}$ in $G$, let $N(\phi)$ be the set of all vertices $u\in V(G)$ such that extending $\phi$ by 
mapping $v$ to $u$ gives a copy of $H_1$ in $G$, and let $n(\phi)=|N(\phi)|$. 
Let $\Phi$ be the set of all copies of $H_1\setminus \{v\}$ in $G$;
then 
$$\sum_{\phi\in \Phi}n(\phi)\ge \vare^{d_1}|G|^{|H_1|}.$$
Let $\Psi$ be the set of all $\phi\in \Phi$ such that $n(\phi)\ge \vare^{d_1+1}|G|$. Since
$$\sum_{\phi\in \Phi\setminus \Psi}n(\phi)\le \sum_{\phi\in \Phi\setminus \Psi} \vare^{d_1+1}|G|\le |G|^{|H_1|-1} \vare^{d_1+1}|G|,$$
it follows that 
$$\sum_{\phi\in \Psi}n(\phi)\ge \vare^{d_1}|G|^{|H_1|}(1-\vare)\ge \vare^{d_1+1}|G|^{|H_1|}.$$
Since $n(\phi)\le |G|$, we deduce that $|\Psi|\ge \vare^{d_1+1}|G|^{|H_1|-1}.$

Let $\phi\in \Psi$.
Thus 
$|N(\phi)|=n(\phi)\ge \vare^{d_1+1}|G|$.
From the choice of $d_2$, either 
there exists $X\subseteq N(\phi)$
with $|X|\ge \vare^{d_2}|N(\phi)|$ such that one of $G[X],\overline{G}[X]$ has edge-density at most $\vare$, or
there are
$\vare^{d_2}|N(\phi)|^{|H_2|}$ copies of $H_2$ in $G[N(\phi)]$.
In the first case, since 
$\vare^{d_2}|N(\phi)|\ge \vare^{d_2}\vare^{d_1+1}|G|\ge \vare^d|G|$, the 
first bullet of (1) holds; so we may assume that 
there are at least 
$$\vare^{d_2}|N(\phi)|^{|H_2|}\ge \vare^{d_2}\vare^{(d_1+1)|H_2|}|G|^{|H_2|}$$ 
copies of $H_2$ in $G[N(\phi)]$, and hence each $\phi\in \Psi$ can be extended to at least $\vare^{d_2}\vare^{(d_1+1)|H_2|}|G|^{|H_2|}$ copies of
$H$. Since $|\Psi|\ge \vare^{d_1+1}|G|^{|H_1|-1}$, there are at least 
$$\vare^{d_1+1}|G|^{|H_1|-1} \vare^{d_2}\vare^{(d_1+1)|H_2|}|G|^{|H_2|} =\vare^{d_1+1+d_2+(d_1+1)|H_2|} |G|^{|H|}= \vare^{d}|G|^{|H|}$$
copies of $H$ in $G$, and hence the second bullet of (1) holds. This proves (1), and hence shows that $H$ is viral, and proves 
\ref{gemsubst2}.~\bbox

Next we will deduce \ref{P4virus}.
The proof uses both \ref{P4thm} and a polynomial bound in the induced graph removal lemma for $P_4$.
The induced graph removal lemma (see \cite{AFKS,conlonfox,Furedi}) says that for each graph $H$ and $0<\vare \leq 1/2$
there exists $\delta>0$ such that every graph $G$ with at most $\delta |G|^{|H|}$ copies of $H$ can be made $H$-free by
adding or deleting at most $\vare |G|^2$ edges. For $H=P_4$, Alon and Fox~\cite{alonfox} proved a polynomial bound:
\begin{thm}\label{P4removal}
There exists $d>0$ such that, if $0<\vare \leq 1/2$, then
every graph $G$ containing at most $\vare^d |G|^{|H|}$ copies of $P_4$ can be made $P_4$-free by
adding or deleting at most $\vare |G|^2$ edges. 
\end{thm}

We deduce \ref{P4virus}, which we restate:
\begin{thm}\label{P4virus2}
Every graph on at most four vertices is viral.
\end{thm}
\Proof Every graph on at most two vertices is viral from the definition. Every graph on three or four vertices,
apart from $P_4$, can be obtained through vertex-substitution from smaller graphs.  So from \ref{gemsubst} it suffices to prove that 
$P_4$ is viral.
Let $d$ be as in \ref{P4removal}, and let $0<\vare\le 1/2$.
We will show that either
\begin{itemize}
\item there exists $X\subseteq V(G)$
with $|X|\ge\vare |G|/4\ge  \vare^3|G|$ such that one of $G[X],\overline{G}[X]$ has at most $\vare\binom{|X|}{2}$ edges; or
\item there are at least $(\vare/4)^{3d}|G|^{4}\ge \vare^{12d}|G|^4$ copies of $P_4$ in $G$.
\end{itemize}
We first dispose of a trivial case, when $\vare |G|/4\le 3$. Then the first bullet holds if $|G|\ge 6$
(because then $G$ or $\overline{G}$ has a triangle), and also if $|G|\le 5$ (because then $\vare |G|/4\le 5/8$ and we can take 
$|X|= 1$). So we may assume that
$\vare |G|/4> 3$.
From the choice of $d$ (with $\vare$ replaced by $(\vare/4)^3$), either $G$ contains at least $(\vare/4)^{3d}|G|^4$
copies of $P_4$ (in which case we are done), or we can obtain a $P_4$-free graph $G'$ with the same vertex set as $G$ by adding or
deleting at most $(\vare/4)^3|G|^2$ edges from $G$. In the latter case, by \ref{P4thm}, there exists $X\subseteq V(G)$
with $|X|\ge \vare |G|/4>3$ such that one of $G'[X],\overline{G'}[X]$ has maximum degree at most $(\vare/4)^2 |G|$. Then
one of $G[X],\overline{G}[X]$ has at most
$$(\vare/4)^3|G|^2 +(\vare/4)^2 |G||X|/2 \leq \frac{3}{8}\vare |X|^2 \leq \vare\binom{|X|}{2}$$
edges (since $|X|\ge 4$). This proves \ref{P4virus2}.~\bbox

\section{Partitioning cographs into R\"odl sets: some counterexamples}

For $\vare>0$, let us say $X\subseteq V(G)$ is {\em $\vare$-restricted} if  one of $G[X],\overline{G}[X]$
has maximum degree at most $\vare|X|$.
There is a strengthening of R\"odl's theorem proved in~\cite{strengthenrodl}:
\begin{thm}\label{strengthenrodl}
For every graph $H$ and every $\vare>0$, there exists $N>0$ such that if $G$ is $H$-free, there is a partition of $V(G)$
into at most $N$ $\vare$-restricted subsets.
\end{thm}
(Note that being $\vare$-restricted involves maximum degree rather than edge-density. The edge-density version is a simple consequence of \ref{rodlthm}.)

There is a corresponding strengthening of the Fox-Sudakov conjecture: perhaps in \ref{strengthenrodl}, $N$ can always
be taken to be a polynomial in $\vare^{-1}$ (depending on $H$).
This seems very intractible, and we have not been able to show it even 
when $H$ is a triangle. But it works when $H=P_4$, as we will show in the next section.

It would be even nicer to get a version of \ref{strengthenrodl} that is strong enough to imply \ref{P4thm} when $H=P_4$, but that 
eludes us. The obvious attempt is false:
\begin{thm}\label{counterex}
For all $\vare$ with $0<\vare<1/2$ such that $\vare^{-1}$ is not an integer, there is a cograph $G$ such that there is no partition of $V(G)$ into at most $1/\vare$ 
$\vare$-restricted sets.
\end{thm}
\Proof
Let $k=\lfloor \vare^{-1}\rfloor$, let $m$ be an integer with $(1-k\vare)m\ge 1$, and let $n$ be some large integer.
Let $G$ be the graph consisting of $k+1$ disjoint cliques $C_0\LL C_k$,
where $|C_0|=n$ and $|C_1|\LL |C_k|=mn$. Suppose that there is a partition of $V(G)$ 
into at most $1/\vare$ 
(and hence at most $k$) $\vare$-restricted sets, and so there is an $\vare$-restricted set $X$
with $|X|\ge mn+n/k$. Let $x_i:=|X\cap C_i|$ for $0\le i\le k$. 
Suppose first that $G[X]$ has maximum degree at most $\vare|X|$. It follows that 
$x_i-1\le \vare|X|$ for $1\le i\le k$, and $x_0\le n$, and summing, 
$$|X|=x_0+x_1+\cdots +x_k\le n+k\vare|X|+k;$$
so 
$$n/(1-\vare k)+n/k\le mn+n/k\le |X|\le (k+n)/(1-k\vare),$$
a contradiction when $n$ is large. Thus, $\overline{G}[X]$ has maximum degree at most $\vare|X|$. Since $|X|>|C_i|$ for $0<i\le k$,
there exists $i\in \{0\LL k\}$ such that $0<x_i\le |X|/2$. Choose $v\in X\cap C_i$; then $v$ has at least $|X|/2>\vare|X|$
non-neighbours in $X$, contradicting that $\overline{G}[X]$ has maximum degree at most $\vare|X|$.
This proves \ref{counterex}.~\bbox

What happens in \ref{counterex} when $\vare^{-1}$ is an integer? Is it true that for every integer $k\ge 1$, every cograph
$G$ can be vertex-partitioned into at most $k$ parts, each $1/k$-restricted? For $k=2$ this is true, and for $k=3$ it is false.
Let us see both those things now.

To show it is true for $k=2$, let us say $X\subseteq V(G)$ is {\em thin} if every component of $G[X]$ has at most
$(|X|+1)/2$ vertices, and {\em thick} if every component of $\overline{G}[X]$ has at most
$(|X|+1)/2$ vertices. So thick and thin sets are both $1/2$-restricted. We will prove:

\begin{thm}\label{2partition} If $G$ is a cograph, there is a partition of $V(G)$ into a thin set and a thick set.
\end{thm}
\Proof
Let $G$ be a cograph. Choose a partition $A,B,C$ of $V(G)$ with $C$ minimal such that 
$|C|> |G|/2$, $A$ is anticomplete to $C$, and $B$ is complete to $C$. (This is possible since we may take $A=B=\emptyset$.)
We may assume that $|C|\ge 2$, and so one of $G[C],\overline{G}[C]$ is not connected, 
by \ref{cograph}. By taking complements if necessary, we may assume that $\overline{G}[C]$ is not connected. Partition
$C$ into two nonempty sets $P,Q$ complete to each other. From the minimality of $C$, it follows that $|P|,|Q|\le  |G|/2$.
In summary, we have a partition of $V(G)$ into four sets $A,B,P,Q$, where $A$ or $B$ may be empty, but $P,Q\ne \emptyset$;
$|A|+|B|<|G|/2$ and $|P|,|Q|\le |G|/2$; $B,P,Q$ are mutually complete, and $A$ is anticomplete to $P\cup Q$. 
(The edges between $A,B$ are unrestricted.)

We may assume that $|Q|\ge |P|$.
Define $m:=\max(0,|Q|-|P|-|B|)$. 
\\
\\
(1) {\em $|A|+m\le 2|Q|$, and $|A|-m\le 2|P|$.}
\\
\\
Suppose first that $m=0$. Then $|Q|\le |P|+|B|$, and we must show that $|A|\le 2|P|$ ($\le 2|Q|$).
But $|P|\ge |Q|-|B|$, so $2|P|\ge |P\cup Q|-|B|\ge |A|$ as required. Now suppose that $m>0$, and so $m=|Q|-|P|-|B|$;
and we must show that $|A|+|Q|-|P|-|B|\le 2|Q|$ and $|A|-(|Q|-|P|-|B|)\le 2|P|$. The first says $|A|-|B|\le |P|+|Q|$, and the 
second that $|A|+|B|\le |P|+|Q|$, and both of these are true. This proves (1).

\bigskip

Consequently, we may choose subsets $P'\subseteq P$ and $Q'\subseteq Q$ with $|P\setminus P'|=\lceil (|A|-m)/2\rceil$ and 
$|Q\setminus Q'|=\lfloor (|A|+m)/2\rfloor$. 
We claim that $X:=A\cup (P\setminus P')\cup (Q\setminus Q')$ is thin and $Y:=B\cup P' \cup Q'$ is thick, and so the theorem holds.
Since $|(P\setminus P')\cup (Q\setminus Q')|=|A|$, and so $|X|=2|A|$, and each of its components has vertex set a subset of 
either $A$ or $(P\setminus P')\cup (Q\setminus Q')$, 
and therefore has
at most $|A|$ vertices, it follows that $X$ is thin. To show that $Y$ is thick, we need:
\\
\\
(2) {\em Each of 
$B, P',Q'$ has cardinality at most $(|Y|+1)/2$.}
\\
\\
Certainly $|B|\le n/2-|A|=|Y|/2$ since $|A|+|B|<n/2$. For the other two inequalities, we have
$|Y|=|P'|+|Q'|+|B|$, and $|P'|=|P|-\lceil (|A|-m)/2\rceil$, and $|Q'|=|Q|-\lfloor (|A|+m)/2\rfloor$, 
so we must show that
$$|P|-\lceil (|A|-m)/2\rceil\le |Q|-\lfloor (|A|+m)/2\rfloor+|B|+1$$
and 
$$|Q|-\lfloor (|A|+m)/2\rfloor\le |P|-\lceil (|A|-m)/2\rceil+|B|+1.$$
These simplify to showing that $|P|+m\le |Q|+|B|+1$, and $|Q|\le |P|+|B|+m$, which both follow from the choice of $m$, and 
since $|Q|\ge |P|$. This proves (2).

\bigskip

From (2), we deduce that $Y$ is thin. This proves \ref{2partition}.~\bbox

Now a counterexample for $k=3$. 
\begin{thm}\label{counterex3}
There is a cograph that admits no vertex-partition into three $1/3$-restricted sets.
\end{thm}
\Proof
Take four disjoint sets $A,B,C,D$ (we will specify their sizes later).
Make $A,B,C$ stable sets complete to each other, and make $D$ a clique anticomplete to $A\cup B\cup C$, forming a graph $G$.
\\
\\
(1) {\em Every $1/3$-restricted subset $X$ of $V(G)$ satisfies either
\begin{itemize}
\item $X\subseteq D$; or
\item $X$ is disjoint from two of $A,B,C$, and its intersection with the third has cardinality at least $2|X\cap D|-3$; or
\item $X$ is disjoint from one of $A,B,C$, and its intersections with the other two and with $D$ have cardinalities
that differ by at most one; or
\item $X\cap D=\emptyset$, and $|X\cap A|, |X\cap B|, |X\cap C|$ differ by at most two; or
\item $|X|\le 6$.
\end{itemize}
}
Suppose first that $\overline{G}$ has maximum degree at most $|X|/3$. Consequently not both
$X\cap (A\cup B\cup C)$ and $X\cap D$ are nonempty, and so we may assume that $X\subseteq A\cup B\cup C$, since otherwise the
first bullet holds. We may assume that $|X\cap A|\ge |X|/3$; but each vertex in $X\cap A$ has at most $|X|/3$ non-neighbours
in $X$, and so $|X\cap A|\le |X|/3+1$; and so the fourth bullet holds.

Now we assume that $G$ has maximum degree at most $|X|/3$. Hence $|X\cap D|\le |X|/3+1$, and so
$|X\cap (A\cup B\cup C)|\ge 2|X|/3-1$. We may assume that $X\cap A\ne \emptyset$; and so $|X\cap(B\cup C)|\le |X|/3$,
and consequently $|X\cap A|\ge |X|/3-1$. If also $X\cap B,X\cap C$ are nonempty, then all three have cardinality
at least $|X|/3-1$ by the same argument, and so vertices in $X\cap A$ have at least $2|X|/3-2>|X|/3$ neighbours in $X$; so $|X|\le 6$
and the fifth bullet holds. So we may assume that $X\cap C=\emptyset$. If also $X\cap B=\emptyset$ then
$X\cap A\ge  2|X|/3-1\ge 2|X\cap  D|-3$ and the second bullet holds, so we assume that $X\cap B\ne\emptyset$. Consequently
$|X\cap A|,|X\cap B|\le |X|/3$, and so $|X\cap D|\ge |X|/3$, and the third bullet holds. This proves (1).

\bigskip

We say that a $1/3$-restricted set $X$ has {\em type} 1--5 depending which bullet of (1) it satisfies.
Now let $n$ be a large integer, and let $|A|=2n$, $|B|=3n$, $|C|=4n$ and $|D|=5n$. Suppose that $V(G)$ can be partitioned into three
$1/3$-restricted sets $X,Y,Z$. If $X$ has type 2, 3, 4 or 5, then $2|X\cap D|\le |X\cap (A\cup B\cup C)|+12$;
so if none of $X,Y,Z$ has type 1, then, summing these three inequalities, we deduce that
$2|D|\le |A|+|B|+|C|+36$, a contradiction since $n$ is large. So we may assume that $Z$ has type 1. If say $Y$ has type
1, 4 or 5, then $|X\cap A|, |X\cap B|,|X\cap C|$ pairwise differ by at least $n-6$, contrary to (1). So $X,Y$ both have types
2 or 3. They cannot both have type 2 since their union includes $A,B,C$, so we assume that $Y$ has type 3. If $X$ has type 2,
then the intersections of two of $A,B,C$ with $X\cup Y$ differ by at most one, a contradiction; so both $X,Y$ have type 3.
But then the sum of two of $|A|,|B|,C|$ should be equal to the third ($\pm 2$), a contradiction. This proves \ref{counterex3}.~\bbox

\section{Partitioning cographs into polynomially many R\"odl sets}

As we just explained, we have not not been able to give a version of \ref{strengthenrodl} that implies \ref{P4thm}.
But at least a ``polynomial'' version of \ref{strengthenrodl} is true when $H=P_4$, because of \ref{P4partn}, which we will now 
prove. The proof breaks into several steps, that follow.

A pair $(P,Q)$ of disjoint subsets of $V(G)$ is {\em pure} if $Q$ is either complete or anticomplete to $P$.
\begin{thm}\label{split}
Let $G$ be a cograph with $|G|\ge 2$, and let $0<\vare\le 1$. Then there is a partition of $V(G)$ into four 
(possibly empty) sets $A_0,A_1,A_2,A_3$, 
with the following properties:
\begin{itemize}
\item $A_0$ is $\vare$-restricted;
\item every two of $A_1,A_2,A_3$ form a pure pair;
\item for $1\le i\le 3$, if $A_i\ne \emptyset$, then there exists $B\subseteq V(G)\setminus A_i$ with $|B|\ge \frac{\vare^2}{4}|G|$ such that 
$(A_i, B)$ is a pure pair.
\end{itemize}
\end{thm}
\Proof For convenience, we say a subset of $V(G)$ is {\em big} if its cardinality is more than $\frac{\vare^2}{4}|G|$, and 
{\em small} otherwise.
Choose a maximal sequence $S_1\LL S_k$ of nonempty, pairwise disjoint, small subsets of $V(G)$, such that
\begin{itemize}
\item for $1\le i\le k$, $S_i$ is complete or anticomplete to $V(G)\setminus (S_1\cupcup S_i)$, and 
\item $|S_1\cupcup S_k|< |G|/2$.
\end{itemize}
Let $A=V(G)\setminus (S_1\cupcup S_k)$. Let $I$ be the set of $i\in \{1\LL k\}$ such that $S_i$ is complete 
to $V(G)\setminus (S_1\cupcup S_i)$, and $J=\{1\LL k\}\setminus I$. Let $P=\bigcup_{i\in I}S_i$
and $Q=\bigcup_{j\in J}S_j$. Thus the sets $S_i\;(i\in I)$ and $A$ are pairwise complete, and the sets 
$S_j\;(j\in J)$ and $A$ are pairwise anticomplete. We may assume that $|A|\ge 2$, since otherwise $|G|\le 2$ and the theorem is true. 
So one of $G[A],\overline{G}[A]$ is not connected, by \ref{cograph}. Choose a partition $B,C$ of $A$ with $B,C$ both nonempty,
such that $B$ is complete or anticomplete to $C$. We may assume that $|B|\ge |C|$, and so $|B|\ge |G|/4$.

From the maximality of the sequence, either:
\begin{itemize}
\item $|B|\le |G|/2$ and $C$ is small; or
\item $C$ is big.
\end{itemize}
In the first case, 
by taking complements if necessary, we may assume that $B,C$ are complete. Let $P'=P\cup C$.
Thus every component of $\overline{G}[P']$ has at most $\frac{\vare^2}{4}|G|$ vertices; so if $|P'|\ge \vare|G|/4$, 
then $P'$ is $\vare$-restricted, and the theorem is satisfied taking $A_0=P'$, $A_1=B$, $A_2=Q$ and $A_3=\emptyset$,
since $B,Q$ are anticomplete.
Note that the third condition of the theorem is satisfied, since $A_1=B$ is complete to the big set $A_0=P'$ 
and $A_2=Q$ is anticomplete to the big set $A_1=B$. 
So we may assume that $|P'|<\vare|G|/4$. Similarly, if $|Q|\ge \vare|G|/4$, then $Q$ is $\vare$-restricted, and
the theorem is satisfied by setting $A_0=Q$, $A_1=B$ and $A_2=P'$, since $A_1$ is anticomplete to the big set $A_0$, and 
$A_2$ is complete
to the big set $A_1$.
But $P'\cup Q=V(G)\setminus B$, and $|B|\le |G|/2$, and so one of $P',Q$ has cardinality at least $|G|/4\ge \vare|G|/4$,
and the theorem holds.

In the second case, since $|B|\ge |G|/4$, \ref{P4thm} implies that there exists $X\subseteq B$ with $|X|\ge \vare|B|$
such that one of $G[X],\overline{G}[X]$ has maximum degree at most $\vare^2|B|$; and by replacing $X$ by a subset, we may assume that
$|X|=\lceil \vare|B|\rceil$. Thus $X$ is $\vare$-restricted.
By taking complements if necessary, we may assume that $G[X]$ has maximum degree at most $\vare|X|$.
Let $Q'=Q\cup X$; then $|Q'|\ge \vare|G|/4$.  If $v\in Q$ then its degree in $G[Q']$ is at most 
$\frac{\vare^2}{4}|G|\le \vare|X|\le \vare|Q'|$; so $Q'$ is $\vare$-restricted, and the theorem is satisfied 
by setting $A_0=Q'$, $A_1=B\setminus X$, $A_2=C$, and $A_3=P$. To see the last, note that  $(A_1,A_2)$ is a pure pair; $A_3$
is complete to both $A_1,A_2$; $A_1$ is complete or anticomplete to the big set $A_2$; $A_2$ is complete or anticomplete to 
the big set $A_1$; and $A_3$ is complete to the big set $A_1$. This proves \ref{split}.~\bbox

Let $X\subseteq V(G)$ with $X\ne\emptyset$. A {\em ribbon attached to $X$} is a sequence $\mathcal{B}=(B_1\LL B_k)$ of pairwise disjoint subsets of $V(G)\setminus X$,
where $k\ge 0$, such that
$B_i$ is complete or anticomplete to $X\cup B_1\cupcup B_{i-1}$ for $1\le i\le k$. Its {\em length} is $k$, and its {\em breadth}
is the minimum of $|B_i|/|X|$ for $1\le i\le k$ (or 1 if $k=0$). We say $X$ is the {\em attachment} of the ribbon. 

We will be concerned with partitions of $V(G)$ into parts, such that for each part $X$, either $X$ is $\vare$-restricted,
or there is a ribbon attached to $X$; and moreover, that for every two parts $X,Y$ that are not $\vare$-restricted, 
$X$ is either complete or anticomplete to $Y$. We must take care that the total number of sets in the partition is not too large,
that the number of beribboned sets is not too large, and that the ribbons are long enough, and have breadth not too small. 
For $\vare>0$, let us say an {\em $(\vare,k)$-beribboning} of a graph $G$ is a partition $\mathcal{P}$ of $V(G)$ (together 
with, implicitly, a choice of ribbons),
such that 
\begin{itemize}
\item for each $X\in \mathcal{P}$, either $X$ is $\vare$-restricted or there is a ribbon of length $k$ attached to $X$; and
\item if $X,Y\in \mathcal{P}$ are different and not $\vare$-restricted, then $(X,Y)$ is a pure pair.
\end{itemize}
(The ribbons attached to distinct members of $\mathcal{P}$ may overlap. We will tidy them up later.) 
The {\em dimensions} of the beribboning are $(m,n)$, where $n=|\mathcal{P}|$ and $m$ is the number of members of $\mathcal{P}$
that are not $\vare$-restricted; and its {\em breadth} is the minimum of the breadth of its ribbons.
From \ref{split}, we see that every cograph admits an $(\vare,1)$-beribboning with dimensions at most $(3,4)$ and breadth at least 
$\vare^2/4$.

\begin{thm}\label{prune}
If $k\ge 0$ is an integer, and $0<\vare\le 1$, and  $G$ is a cograph that admits an $(\vare,k)$-beribboning with dimensions at most $(m,n)$ and breadth at least $\beta$, then it also 
admits  an $(\vare,k)$-beribboning with dimensions at most $(\vare^{-2},n)$ and breadth at least $\beta$.
\end{thm}
\Proof We will prove that if $G$ admits an $(\vare,k)$-beribboning $\mathcal{P}$ with dimensions $(m,n)$ and breadth $\beta$, where
$m\ge \vare^{-2}$, then it also admits one with dimensions at most $(m-1,n)$ and breadth at least $\beta$.
Let $t=\lceil \vare^{-1}\rceil$. Since $m\ge \vare^{-2}>(t-1)^2$, and every cograph with more than $(t-1)^2$
vertices has a clique or stable set of size $t$, we may choose distinct 
$X_1\LL X_{t}\in \mathcal{P}$, not $\vare$-restricted, such that either $X_1\LL X_{t}$ are pairwise anticomplete 
or $X_1\LL X_{t}$ are pairwise complete. We may assume that $|X_t|\le |X_1|\LL |X_{t-1}|$. Choose $Y_i\subseteq X_i$
with $|Y_i|=|X_t|$ for $1\le i\le t-1$. Then $Y_1\cupcup Y_{t-1}\cup X_t$ is $\vare$-restricted (because $t\ge \vare^{-1}$).
Let $\mathcal{P}'$ be the partition obtained from $\mathcal{P}$ by replacing the sets $X_1\LL X_t$ by the sets
$$X_1\setminus Y_1\LL X_{t-1}\setminus Y_{t-1}, Y_1\cupcup Y_{t-1}\cup X_t.$$
Thus $\mathcal{P}'$ has the same number of members as $\mathcal{P}$, but at least one more of them is $\vare$-restricted. Moreover,
each of the sets $X_i\setminus Y_i$ has a ribbon attached of length $k$ and breadth at least $\beta$. So $\mathcal{P}'$
is  an $(\vare,k)$-beribboning with dimensions at most $(m-1,n)$ and breadth $\ge \beta$. By repeating, this proves \ref{prune}.~\bbox

\begin{thm}\label{growtree}
If $G$ is a cograph, and $0<\vare\le 1$, and $k\ge 0$ is an integer, then $G$ admits an $(\vare,k)$-beribboning with dimensions 
at most $(\vare^{-2},1+3k\vare^{-2})$ and breadth at least $\vare^2/4$.
\end{thm}
\Proof
We proceed by induction on $k$. For $k=0$, the result is trivial.
Inductively, we assume that $k\ge 1$, and 
$G$ admits an $(\vare,k-1)$-beribboning $\mathcal{P}$ with dimensions at most $(\vare^{-2},1+3(k-1)\vare^{-2})$ and breadth 
at least $\vare^2/4$. Let $X_1\LL X_m$ be the members of $\mathcal{P}$ that are not $\vare$-restricted.
For $1\le i\le m$, since $|X_i|\ge 2$, there is, by \ref{split},
a partition of $X_i$ into four
(possibly empty) sets $A_{i0},A_{i1},A_{i2},A_{i3}$,
with the following properties:
\begin{itemize}
\item $A_{i0}$ is $\vare$-restricted;
\item every two of $A_{i1},A_{i2},A_{i3}$ form a pure pair;
\item for $1\le j\le 3$, if $A_{ij}\ne \emptyset$, then there exists $B\subseteq X_i\setminus A_{ij}$ such that
$|B|\ge \frac{\vare^2}{4}|X_i|$ and $(A_{ij}, B)$ is a pure pair.
\end{itemize}
Let $\mathcal{P}'$ be obtained from $\mathcal{P}$ by replacing $X_i$ by $A_{i0},A_{i1},A_{i2},A_{i3}$ for $1\le i\le m$.
We claim that this is a $(\vare,k)$-beribboning with dimensions
at most $(3m,1+3k\vare^{-2})$ and breadth at least $\vare^2/4$. Certainly 
$$|\mathcal{P}'|\le |\mathcal{P}|+3m\le 1+3k\vare^{-2}$$
and the number of members of $\mathcal{P}'$ that are not $\vare$-restricted is at most $3m$. We need to check the ribbons.
Let $1\le i\le m$. There is a $(\vare,k-1)$-ribbon $(B_1\LL B_{k-1})$ attached to $X_i$ with breadth at least $\vare^2/4$.
Let $1\le j\le 3$ with $A_{ij}\ne \emptyset$. From the third bullet above, there exists 
a subset $B\subseteq X_i\setminus A_{ij}$ such that 
$|B|\ge \frac{\vare^2}{4}|X_i|$ and $(A_{ij}, B)$ is a pure pair. But then $(B,B_1\LL B_{k-1})$
is a ribbon attached to $A_{ij}$ of breadth at least $\vare^2/4$ and length $k$. This proves our claim that
$\mathcal{P}'$ is a $(\vare,k)$-beribboning with dimensions
at most $(3m,1+3k\vare^{-2})$ and breadth at least $\vare^2/4$, and then an application of \ref{prune} gives the result.
This proves \ref{growtree}.~\bbox

Let us say a ribbon $(B_1\LL B_k)$ attached at $X$ is {\em pure} if either all the sets $X,B_1\LL B_k$ are pairwise complete,
or all the sets $X,B_1\LL B_k$ are pairwise anticomplete. An $(\vare,k)$-beribboning is {\em pure} if all the ribbons 
it uses are pure.
\begin{thm}\label{pureribbon}
If $G$ is a cograph, and $0<\vare\le 1/2$, then $G$ admits a pure $(\vare,\lceil\vare^{-1}\rceil)$-beribboning 
with dimensions
at most $(\vare^{-2},10\vare^{-3})$ and breadth at least $\vare^2/4$.
\end{thm}
\Proof
Taking $k=2\lceil \vare^{-1}\rceil\le 3\vare^{-1}$, we deduce from \ref{growtree} that $G$ admits an
$(\vare,k)$-beribboning with dimensions
at most $(\vare^{-2},1+3k\vare^{-2})\le (\vare^{-2},10\vare^{-3})$ and breadth at least $\vare^2/4$. Let $(B_1\LL B_k)$ be a ribbon attached to some 
$X\subseteq V(G)$. Let $I$ be the set of $i\in \{1\LL k\}$ such that $B_i$ is complete to $X\cup B_1\cupcup B_{i-1}$,
and $J=\{1\LL k\}\setminus I$. Then both of $(B_i:i\in I)$, $(B_j:j\in J)$
are pure ribbons attached to $X$, and one of them has length at least $k/2=\lceil \vare^{-1}\rceil$. Hence, for 
 each $X\in \mathcal{P}$ that is not $\vare$-restricted, there is a 
pure ribbon of length at least $\vare^{-1}$ and breadth at least $\vare^2/4$ attached to  $X$.
This proves \ref{pureribbon}.~\bbox

The purpose of the ribbons is: suppose we have a pure ribbon $(B_1\LL B_k)$ attached to a set $X$, and its length is at 
least $\vare^{-1}$, and $B_1\LL B_k$ all have the same size. Then we can partition $X\cup B_1\cupcup B_k$ into $\vare$-restricted 
sets, as follows.
We can partition almost all of $X$ into a few $\vare$-restricted subsets, greedily, in such a way that the remainder, $Y$ say, 
has size at most $|B_1|$;
and then $Y\cup B_1\cupcup B_k$ is also $\vare$-restricted. 
But to use this method, we first need
to tidy up the ribbons. We want to arrange that:
\begin{itemize}
\item the ribbons are vertex-disjoint from one another;
\item for each set of the partition $\mathcal{P}$, at most half its vertices belong to ribbons; and 
\item for each ribbon, all its members have the same size.
\end{itemize}
We call these the {\em prettification conditions}.
Let us say a $(\vare,k)$-beribboning is {\em prettified} if it is pure and satisfies the three conditions above.
All these things will be accomplished by replacing the sets of the ribbons by subsets of themselves. This will reduce the 
breadth, so we must be careful that the breadth does not get too small. In particular, if the sets of a ribbon are already 
very small, we may not be able to shrink them by the required factors, and we must treat such ribbons differently. But
in this case, the corresponding attachment is also very small, a polynomial in $\vare^{-1}$, and we can easily partition it
into a few $\vare$-restricted sets, and need not use the ribbon at all.

Let us see the last statement above. By repeatedly applying \ref{P4thm}, we evidently have:
\begin{thm}\label{smallset}
If $G$ is a cograph, and $X\subseteq V(G)$, and $0\le \vare\le 1$, and $t\ge 0$ is an integer, then there are $t$ pairwise disjoint $\vare$-restricted 
subsets of $X$ with union $X\setminus Y$ say, such that $|Y|\le (1-\vare)^t|G|\le e^{-\vare t}|X|$.
\end{thm}
So in particular, there is a partition of $X$ into $t$ $\vare$-restricted sets, if
we take $t$ so large that $e^{\vare t}>|X|$.

\begin{thm}\label{prettify}
For $0<\vare\le 1/2$, if $G$ is a cograph, then 
$G$ admits a prettified $(\vare,\lceil \vare^{-1}\rceil)$-beribboning $\mathcal{P}$ with dimensions at most $(\vare^{-2},21\vare^{-4})$ and 
breadth  at least $\vare^4/32$.
\end{thm}
\Proof Let $k=\lceil \vare^{-1}\rceil$. By \ref{pureribbon}, $G$ admits a pure $(\vare,k)$-beribboning $\mathcal{P}_0$
with dimensions at most 
$(\vare^{-2},10\vare^{-3})$ and breadth at least $\vare^2/4$.
Let $q=16$.
If $X\in \mathcal{P}_0$ and $|X|<\vare^{-q}$, then by \ref{smallset} there is a partition of $X$ into at most 
$q\vare^{-1}\log(\vare^{-1})\le q\vare^{-2}$ $\vare$-restricted
sets. For each $X\in \mathcal{P}_0$ that is not $\vare$-restricted and has cardinality less than $\vare^{-q}$,
let us replace $X$ in $\mathcal{P}_0$ 
by the sets of the corresponding partition. We obtain a 
pure $(\vare,k)$-beribboning $\mathcal{P}$ of $G$ with dimensions at most 
$(\vare^{-2},10\vare^{-3}+q\vare^{-4})\le (\vare^{-2},(q+5)\vare^{-4})$ and breadth at 
least $\vare^2/4$,
such that each $X\in \mathcal{P}$ that is not $\vare$-restricted has cardinality at least $\vare^{-q}$. 
Now we will prettify the ribbons.

Let $X_1\LL X_m$ be the members of $\mathcal{P}$ that are not $\vare$-restricted, and for $1\le i\le m$, let $(B_{i1}\LL B_{ik})$
be a pure ribbon attached to $X_i$ with breadth at least $\vare^2/4$. 
\\
\\
(1) {\em For $1\le i\le m$ and $1\le j\le k$ there exists $C_{ij}\subseteq B_{ij}$, with $|C_{ij}|\ge \lfloor |B_{ij}|/m\rfloor$
such that all the sets $C_{ij}\; (1\le i\le m, 1\le j\le k)$ are pairwise disjoint.}
\\
\\
For each $i,j$, let $\mathcal{Q}_{ij}$ be a partition of $B_{ij}$ into $\lceil |B_{ij}|/m\rceil$ sets, 
of which $\lfloor |B_{ij}|/m\rfloor$ have size $m$, and possibly one has size less than $m$.
For each $Q\in \mathcal{Q}_{ij}$ take a new vertex $u_{ij}^Q$, and let $U$ be the set of all these new vertices; that is,
$$U=\{u_{ij}^Q: 1\le i\le m, 1\le j\le n, Q\in \mathcal{Q}_{ij}\}.$$
Let $H$ be the graph with bipartition $(U,V(G))$
in which $u_{ij}^Q$ is adjacent to $v\in V(G)$ if $v\in Q$. Then $H$ has maximum degree at most $m$, and so 
its edge-set can be partitioned into $m$ matchings. Let $M$ be one of these matchings.
For each vertex $h\in V(H)$ with degree $m$ in $H$, there is an edge of $M$
incident with $h$. Consequently, for $1\le i\le m$ and $1\le j\le k$, there are at least $\lfloor |B_{ij}|/m\rfloor$
edges in $M$ with ends in $\mathcal{Q}_{ij}$.
Let $C_{ij}$ be the set of vertices of $G$ that are joined by an edge of $M$ to a vertex in $\mathcal{Q}_{ij}$.
This proves (1).

\bigskip

For $1\le i\le m$, since $|X_i|\ge \vare^{-q}$, it follows that each $|B_{ij}|\ge \vare^{2-q}/4\ge \vare^{-2}\ge m$, and hence
$$|C_{ij}|\ge \lfloor|B_{ij}|/m\rfloor\ge |B_{ij}|/(2m)\ge (\vare^4/8)|X_i|.$$
Consequently
 $(C_{i1}\LL C_{ik})$ is a ribbon attached to $X_i$ with breadth at least $\vare^4/8$.
These ribbons are pairwise vertex-disjoint, so we have satisfied the first prettification condition.

For the second condition, for each $i,j$, choose $D_{ij}\subseteq C_{ij}$ such that for each $Y\in \mathcal{P}$,
$|D_{ij}\cap Y|=\lfloor |C_{ij}\cap Y|/2\rfloor$. Then 
$$|D_{ij}|\ge |C_{ij}|/2-|\mathcal{P}|/2\ge |C_{ij}|/2-(q+5)\vare^{-4}/2\ge |C_{ij}|/4\ge (\vare^4/32)|X_i|.$$
since 
$$|C_{ij}|\ge (\vare^4/8)|X_i|\ge \vare^{4-q}/8\ge    2(q+5)\vare^{-4}.$$
Then for each $Y\in \mathcal{P}$, at most half the vertices of $Y$ belong to the union of the sets $D_{ij}$; 
so the ribbons $(D_{i1}\LL D_{ik})$ have breadth at least $\vare^4/32$ and satisfy the first and second prettification conditions. 
Finally,
for the third, for each $i,j$ choose $E_{ij}\subseteq D_{ij}$ of cardinality $\lceil (\vare^4/32)|X_i|\rceil$; then
the ribbons $(E_{i1}\LL E_{ik})$ for $1\le i\le m$ satisfy all the prettification conditions, and have breadth at least
$\vare^4/32$. This proves \ref{prettify}.~\bbox

Finally, let us deduce \ref{P4partn}, which we restate, with $\vare$ replaced by $\vare/2$ for convenience:
\begin{thm}\label{P4partn2}
If $G$ is a cograph, then for every $\vare$ with $0<\vare\le 1/2$,
there is a partition of $V(G)$ into at most $30\vare^{-4}$ sets, all $2\vare$-restricted.
\end{thm}
\Proof
From \ref{prettify}, 
$G$ admits a prettified $(\vare,\lceil \vare^{-1}\rceil)$-beribboning $\mathcal{P}$ with dimensions at most $(\vare^{-2},21\vare^{-4})$ and
breadth  at least $\vare^4/32$. Let $X_1\LL X_m$ be the members of $\mathcal{P}$ that are not $\vare$-restricted, 
and for $1\le i\le m$ let $(B_{i1}\LL B_{ik})$
be a pure ribbon attached to $X_i$ with breadth at least $\vare^4/32$, satisfying the prettification conditions.
Let $F$ be the union of the sets $B_{ij}$ for $1\le i\le m$ and $1\le j\le k$. Thus $|F\cap Y|\le |Y|/2$
for each $Y\in \mathcal{P}$.  Let $Z=F\cup X_1\cupcup X_m$. We will partition $Z$ and $V(G)\setminus Z$ into
$2\vare$-restricted sets.

Let $t= \lceil 8\vare^{-2}\rceil\le 9\vare^{-2}$. By \ref{smallset}, for $1\le i\le m$, since 
$$|X_i\setminus F|\le |X_i|\le (32\vare^{-4})|B_{i1}|,$$ 
we may choose $t$
pairwise disjoint $\vare$-restricted subsets $A_{i1}\LL A_{it}$ of $X_i\setminus F$, with union 
$(X_i\setminus F)\setminus Y_i$ say, 
such that $|Y_i|\le |B_{i1}|$. But then the sets 
\begin{align*}
&A_{ij}\; (1\le i\le m, 1\le j\le t)\\
&Y_i\cup B_{i1}\cupcup B_{ik}\; (1\le i\le m)
\end{align*}
are all $\vare$-restricted, and pairwise disjoint, and have union 
$Z$. 

For each $\vare$-restricted set $Y\in \mathcal{P}$, the set $Y\setminus F$ is
$2\vare$-restricted, since $|Y\setminus F|\ge |Y|/2$; and the union of all these sets $Y\setminus F$ (where $Y\in \mathcal{P}$
is not $\vare$-restricted) is $V(G)\setminus Z$.
So altogether we have found a partition of $V(G)$ into at most 
$$|\mathcal{P}|+m t\le 21\vare^{-4} + \vare^{-2}(9\vare^{-2})\le 30\vare^{-4}$$
$2\vare$-restricted sets. This proves \ref{P4partn2}.~\bbox

\section*{Acknowledgement}
Thanks to Maria Chudnovsky and Sophie Spirkl for stimulating discussions on the topic of this paper.

\end{document}